\documentclass[preprint,12pt]{elsarticle}

\usepackage{amsthm,amssymb,amsmath}
\usepackage{subcaption}
\usepackage{caption}
\usepackage{graphicx}
\usepackage{stmaryrd}
\usepackage{mathrsfs}
\usepackage{bm}
\usepackage{placeins}
\usepackage{textgreek}
\usepackage{xurl}
\usepackage{hyperref}
\usepackage{float}
\usepackage{placeins}
\usepackage{siunitx}
\usepackage{geometry}
\usepackage{cleveref}
\usepackage{multirow}
\usepackage{mathtools}
\usepackage{amsthm}
\usepackage{amssymb}
\usepackage{subcaption}
\usepackage{natbib}
\usepackage{xcolor}
\usepackage{mathpazo}
\usepackage{float}
\usepackage{adjustbox}
\usepackage{enumitem}
\usepackage{lipsum}
\usepackage[most]{tcolorbox}
\usepackage{orcidlink}
\numberwithin{equation}{section}

\definecolor{lavender}{RGB}{230, 230, 250}
\definecolor{peach}{RGB}{255, 229, 180}
\definecolor{mintcream}{RGB}{245, 255, 250}
\hypersetup{
	colorlinks=true,
	linkcolor=cyan,
	citecolor=cyan,
	urlcolor=cyan}
\usepackage{cleveref}
\usepackage{color}
\usepackage{titlesec}
\usepackage[T1]{fontenc}
\usepackage{stmaryrd}
\usepackage{lmodern}

\newtheorem{theorem}{Theorem}[section]

\newtheorem{lemma}{Lemma}[section]
\theoremstyle{definition}
\newtheorem{example}{Example}[section]

\usepackage{amsthm}
\usepackage{titlesec}

\titleformat{\subsection}
{\normalfont\normalfont\bfseries}
{\thesubsection}{1em}{}
\newtheoremstyle{myremark}
{6pt}   
{6pt}   
{\normalfont} 
{}      
{\fontseries{sb}\selectfont} 
{.}     
{0.5em} 
{\thmname{#1}\ \thmnumber{#2}} 

\theoremstyle{myremark}

\journal{Elsevier journal}

\begin{document}
	
	\begin{frontmatter}
		
		
		
		
		\title{\textbf{Entropy bounded compact stencil Runge-Kutta discontinuous Galerkin method for Euler equations}}
		
		
		\author[]{Deepak Kumar\corref{cor1} \orcidlink{https://orcid.org/0009-0008-3829-7985}}
		
	\author[]{Pratima Rai~\orcidlink{https://orcid.org/0000-0003-3872-1148}}
		
	\author[]{Arvind Patel \orcidlink{https://orcid.org/0000-00001-7345-6859} }

		\address{Department of Mathematics, University of Delhi, New Delhi 110007, India}

	\cortext[cor1]{Corresponding author. Email: dkuanr.phd2025@maths.du.ac.in}
		
		
		\setlength{\parskip}{3pt}
		
		\begin{abstract}
		 In this paper we introduce a compact stencil Runge-Kutta discontinuous Galerkin method for the Euler equations of compressible gas dynamics, which preserves a minimum entropy principle similar to that introduced by Tadmor for entropy solution of the Euler equations \cite{Tadmor}. Compact stencil RKDG, introduced by Chen et al. in \cite{compact rkdg}, is a variant of standard RKDG method which is known for its smaller stencil size and parallel efficiency. The main challenge is that the cRKDG method has to be expressed in Butcher form and not every Runge-Kutta method of this form is a convex combination of forward-Euler steps. Therefore, in order to establish this property we leverage the bound preserving framework of the cRKDG schemes \cite{bound preserving crkdg}, an entropy bounded limiter \cite{Lv Yu} and the classical positivity preserving limiter. In the last section numerical experiments are provided to demonstrate the robustness of the method.
		\end{abstract}
		
		\clearpage

		\begin{highlights}
			\item A minimum entropy principle for the Lax-Friedrichs scheme with a definite CFL condition is proven
	   \item Entropy bounded compact stencil RKDG is constructed
	   \item Bound preserving framework of cRKDG is employed for this purpose
		\end{highlights}
		
		
	\end{frontmatter}
	
		
		
		\section{Introduction}
		\label{sec1}
	Hyperbolic conservation laws are a class of partial differential equations arising in several areas of scientific interest such as shallow water waves, astrophysics, atmospheric flows, magnetohydrodynamics etc. An interesting feature of this type of differential equation is that its solution can become discontinuous in finite time even if the associated initial condition is smooth. Therefore, it becomes significantly challenging to simulate their solutions using high order numerical schemes which exhibit Gibbs oscillations near discontinuities. \newline
	Discontinuous Galerkin method is a widely used Finite element technique known for its high order accuracy, compactness, parallel-scalability, ability to handle complex geometries, $hp$-adaptivity etc. If Runge-Kutta time integration is followed by discontinuous Galerkin spatial discretization then the resulting method is a type of standard method of lines and is known as the Runge-Kutta discontinuous Galerkin method which is often abbreviated as RKDG \cite{cockburn1998, cockburn1990, cockburn1989a}.
	For decades it has been one of the widely used numerical technique for solving hyperbolic conservation laws. However, the method often needs limiting or extra stabilization mechanism to restrain spurious oscillations and instabilities. Some well established limiting techniques includes the use of slope limiters \cite{cockburn1989b}, WENO limiters \cite{Weno limiters}, artificial viscosity \cite{chan} etc. A more efficient approach to mitigate spurious oscillations without relying on any parameter based identification of troubled cells is to damp the higher frequency modes of the polynomial solution \cite{lu2021, lu2022, OEDG}. \\
	Recently, Chen et al.\cite{compact rkdg} proposed compact stencil Runge-Kutta discontinuous Galerkin method for hyperbolic conservation laws. It departs from the standard method of lines and introduces a local derivative operator in the inner stages of the RK method without destroying the accuracy and local conservation of the original RKDG method. A similar philosophy of hybridizing different spatial operators at different stages of time integration is further explored in stage dependent RKDG \cite{sdrkdg} and our recent work \cite{deepak}. \\
	Tadmor proved that the entropy solution of Euler equations satisfy a minimum entropy principle \cite{Tadmor}. For higher order DG schemes this property was established in \cite{Zhang entropy, Lv Yu}. However, their work rely on structure of SSP-RK to extend the entropy boundedness of the DG scheme with forward-Euler update to a DG scheme with high order explicit RK. Therefore, the similar arguments cannot be applied to the cRKDG method which could only be expressed in Butcher form. In this paper, we derive the entropy bounded principle for the cRKDG schemes by leveraging the bound-preserving framework of the cRKDG schemes. We first derive a minimum-entropy principle for the Lax-Friedrichs scheme which is then used to deduce a similar entropy principle for the forward-Euler update with local derivative operator. Finally, we express each stage of the cRKDG method as a convex combination of the forward-Euler updates involving the DG operator and the local derivative operator. The quasi-concavity of specific thermodynamic entropy yields the desired result. In the last section, we have provided numerical experiments to demonstrate the accuracy and robustness of entropy-bounded cRKDG method.
		\section{Scheme Formulation}
We consider the conservation law of the form 
\begin{align}
	\mathbf{u}_{t} + \nabla \cdot \mathbf{f}(\mathbf{u}) = 0 \qquad \text{in} \quad \Omega, 
	\label{equation:conservation law}
\end{align}
where $\Omega$ is an open subset of $\mathbb{R}^d$.  
It is a well known fact that the solution to these problems can become discontinuous in finite time. In this scenario the notion of classical solution is no longer valid and one has to move towards the associated variational form and it's solution which is regarded as the weak solution. However, there can be infinitely many weak solutions and additional admissibility criterions are needed to single out the physically relevant solution which is also known as the "entropy solution". One such criterion is the entropy inequality which is always satisfied by the physically admissible solution. In other words, if $\mathbf{u}$ is an entropy solution of \eqref{equation:conservation law}, then for any entropy function $ \mathcal{U}$, and the corresponding entropy flux function $\mathcal{F} $ the following inequality holds in the sense of distributions:
\begin{align}
	\frac{\partial \mathcal{U}(\mathbf{u})}{\partial t} + \nabla \cdot \mathcal{F}(\mathbf{u}) \leq 0 ,
	\label{equation:entropy inequality}
\end{align}
Note that the compatibility condition $ \frac{\partial \mathcal{F}}{\partial \mathbf{u}} = (\mathcal{U}^{\prime} )^{T}  \frac{\partial \mathcal{\mathbf{f}}}{\partial \mathbf{u}}$ along with \eqref{equation:entropy inequality} has to be satisfied. A necessary and sufficient condition for any convex function $\mathcal{U}$ to be an entropy function is the symmetry of the matrix $ \mathcal{U}^{\prime \prime}(\mathbf{u}) \mathbf{f}^{\prime}(\mathbf{u}).$
A particular case of conservation laws \eqref{equation:conservation law} that is of much interest is the Euler equations of compressible gas dynamics where the vector of conserved variables and flux function are respectively given by
\begin{align}
	\mathbf{u} = \left(\rho, \rho \mathbf{v}^T, \rho e\right)^{T} \quad \text{and} \quad \mathbf{f(\mathbf{u})} = (\rho \mathbf{v}^T, \rho \mathbf{v} \otimes \mathbf{v} + p I, \mathbf{v}^T(\rho e + p) )^{T}.
\end{align}
Here $\rho$ denotes the density, $\mathbf{v}$ is velocity, $e$ is specific total energy and $p$ is the pressure. The equation of state for the ideal gas gives the relation between pressure and internal energy. i.e.
\begin{align*}
	p = (\gamma - 1) \left(\rho e - \frac{\rho |\mathbf{v}|^2}{2} \right),
\end{align*}
where $\gamma$ is the ratio of specific heat capacities and the speed of sound is defined as $ c = \sqrt{\frac{\gamma p}{\rho}} $. A well known entropy-entropy flux pair for the Euler equations is $ \left(\mathcal{U}, \mathcal{F}\right) = \left(-\rho s, -\rho \mathbf{v} s\right)$, where $s$ is the specific thermodynamic entropy given by
\begin{align}
	s = \ln(p) - \gamma \ln(\rho) + \delta. 
	\label{equation:physical entropy}
\end{align}
Note that for the physical entropy to be well defined, the pressure and density are required to be positive. For Euler equations of gas dynamics the solution is sought in a convex subset of $\mathbb{R}^m$ so that the pressure and density are always positive. i.e.
\begin{align}
	G_0 = \{  \mathbf{u} = \begin{pmatrix}
		\rho \\
		\rho \textbf{v} \\
		\rho e
	\end{pmatrix} : \rho > 0 \quad \text{and} \quad p > 0  \}
\end{align}
The convexity of this set follows from the concavity of pressure and the fact that convex combination of positive densities is positive. The high order schemes on the other hand are not positivity preserving vacuously and during the numerical simulations the density and pressure can often become negative. However, the element averaged solution at the succeeding time level can still be proven positivity preserving if the solution at the current time level lies in the set $ G_0$. Of course the time evolution method chosen should be a convex combination of the forward-Euler steps and a constraint on the time step size is required. In \cite{Tadmor} the author proved that the entropy solution of the Euler equations of gas dynamics satisfy a minimum entropy principle. In addition, the author also proved that the three point finite-difference/volume approximations of these equations satisfy a similar entropy principle under a suitable constraint on the time step size. The framework was generalized for the high order finite volume schemes and the discontinuous Galerkin method in \cite{Zhang entropy}, by Zhang et al. They have shown that the specific entropy is quasi-concave and the element averaged solutions of these high order schemes lie in the following convex set
\begin{align}
	G = \{ \mathbf{u} = \begin{pmatrix}
		\rho \\
		\rho \textbf{v} \\
		\rho e
	\end{pmatrix} : \rho > 0, p > 0 \quad \text{and} \quad  s(\mathbf{u}) \ge  \min_x s_0(\mathbf{x})\},
\end{align}
where $s_0$ denotes the entropy of the initial condition. We would like to briefly mention their approach. For DG scheme or high order FV scheme they have shown that the forward Euler update of the element averaged solution is bounded below by the minimum entropy of the current solution at the quadrature points of the element and it's values at the interfaces. A scaling type limiter similar to the context of bound preserving techniques is then used to enforce the entropy boundedness of the polynomial solution. However, the scaling parameter for every element is determined by the Newton-Raphson method. A more efficient approach to avoid the expensive iterative procedure and to construct an entropy-bounded DG scheme for the Euler equations was proposed in \cite{Lv Yu}. The generalization to DG scheme with high order time integration follows from the structure of SSP-RK, which can be expressed as convex combination of the forward Euler updates. The scaling type limiter being applied after each stage of the Runge-Kutta method. For the compact RKDG method this generalization is not straightforward because the method has to be expressed in Butcher form to preserve local conservation. Therefore, the earlier methodology cannot be applied to the cRKDG scheme without careful consideration of the entropy boundedness at every stage.
\section{Notations and preliminaries}
Let us assume that the domain $\Omega$ is triangulable and one can partition it into a set of non-overlapping elements. We consider any such partition of $\Omega$ and denote it by $\mathcal{T}_h$ where  $h = \max_{K \in \mathcal{T}_h}h_K$ is the mesh parameter. In addition, we also define the finite element space 
\begin{align}
	\mathcal{V}_h^k = \{v \in L^2(\Omega): v|_K \in [P^k(K)]^m, \forall K \in \mathcal{T}_h\},
\end{align}
where $P^k(K)$ is the space of polynomial of degree at most $k$ on the element $K.$ Let $n_{\ell, K}$ be the unit outward normal to the edge $e_{\ell, K}$ of $K.$ For any member $v$ of $\mathcal{V}_h^k$, we will use the notations $v^{-} \, \text{and} \, v^{+} $ to denote the interior and exterior traces of $v$ along any edge $e_{\ell, K}.$ 
The discontinuous Galerkin method for solving \eqref{equation:conservation law} in semidiscrete form is then given by:- find $\mathbf{u}_h  \in \mathcal{V}_h^k$ such that for any test function $v \in \mathcal{V}_h^k$, the following relation holds
\begin{align}
	\int_K\left(\mathbf{u}_h\right)_t: v_h \, dx -\int_K\mathbf{f}(\mathbf{u} _h) : \nabla v_h\, d x+\sum_{\ell = 1}^{l(K)} \int_{e_{\ell, K}} \widehat{\mathbf{f} \cdot  n_{\ell, K}}(\mathbf{u}_h^{-}, \mathbf{u}_h^{+}) \cdot  \ v_h \, d s = 0, 
	\label{equation:semiDgscheme}
\end{align}
where ":" is the Forbenius inner product and $\widehat{\mathbf{f} \cdot  n_{\ell, K}}$ is the numerical flux, which solves a Riemann problem exactly or approximately at each interface. A well known example of numerical flux is local Lax-Friedrichs flux which is given by 
\begin{align}
 \widehat{\mathbf{f} \cdot  n_{\ell, K}}(\mathbf{u}_h^{-}, \mathbf{u}_h^{+}) = 
\frac{ \mathbf{f}(\mathbf{u}_h^{-}) + \mathbf{f}(\mathbf{u}_h^{+})}{2} \cdot  n_{\ell, K} - \alpha_{\ell, K}\frac{(\mathbf{u}_h^{+} - \mathbf{u}_h^{-})}{2} \quad \text{where} \quad \alpha_{\ell, K} = \max | eig(\mathbf{f}^{\prime} \cdot n_{\ell, K})| .
\nonumber
\end{align}
We will use Lax-Friedrichs numerical flux in this paper and dissipation parameter for this will be denoted by $\alpha_0 = \max_{x, \ell, K} \alpha_{\ell, K}.$ Next, we define the DG operator $\nabla^{DG}\cdot \mathbf{f}: \mathcal{V}_h^k \rightarrow \mathcal{V}_h^k $ as follows:
\begin{align}
   \int_K \nabla^{DG} \cdot \mathbf{f}(\mathbf{u}_h)\cdot v_h\, dx 
	= -\int_K \mathbf{f}(\mathbf{u}_h) : \nabla v_h\, dx +\sum_{\ell = 1}^{l(K)} \int_{e_{\ell, K}} \widehat{\mathbf{f} \cdot  n_{\ell, K}}(\mathbf{u}_h^{-}, \mathbf{u}_h^{+}) \cdot  \ v_h \, d s, 
\end{align}
 and the local derivative operator $\nabla^{loc}\cdot \mathbf{f}: \mathcal{V}_h^k \rightarrow \mathcal{V}_h^k $ is defined by
\begin{align}
\int_K \nabla^{loc} \cdot \mathbf{f}(\mathbf{u}_h)\cdot v_h\, dx 
	= -\int_K \mathbf{f}(\mathbf{u}_h) : \nabla v_h\, dx  +\sum_{\ell = 1}^{l(K)}\int_{e_{\ell, K}} \left(\mathbf{f}(\mathbf{u}_h^{-}) \cdot n_{\ell, K}\right)  \cdot  \ v_h \, d s.
\end{align}
For convenience we will use the notations $  \mathscr{F} (\mathbf{u}_h) = -\nabla^{DG} \cdot \mathbf{f}(\mathbf{u}_h) $ and $  \mathcal{G} (\mathbf{u}_h) = -\nabla^{loc} \cdot \mathbf{f}(\mathbf{u}_h).$ A Butcher form of RK can be used to further discretize \eqref{equation:semiDgscheme} in time:
 \[
 \quad \begin{array}{c|c}
	c & A\\
	\hline
	& b
\end{array}, \quad A = \left(a_{ij}\right)_{s\times s}, \quad  b = \left(b_1, b_2,\cdots b_s\right), \quad  c = \left(c_1, c_2,\cdots c_s\right)^T, \]
where $A $ denotes the coefficient matrix which is strict lower triangular and $b$ denotes the weights. If $\Delta t$ is the time step size then the fully discrete RKDG method is given by
	\begin{equation}
	\begin{aligned}
		& \mathbf{u}_h^{(i)}= \mathbf{u}_h^n + \Delta t \sum_{j=1}^{i-1} a_{i j} \mathscr{F}(\mathbf{u}_{h}^{(j)} ) , \quad i = 1,2\hdots ,s,  \\
		\label{equation:RKDG fully discrete}
		&\mathbf{u}_h^{n+1} = \mathbf{u}_h^n + \Delta t \sum_{i = 1}^{s} b_\ell \mathscr{F}(\mathbf{u}_{h}^{(i)} ),  \\  
	\end{aligned}
\end{equation}
The cRKDG method can be obtained from equation \eqref{equation:RKDG fully discrete} by replacing the DG operator in the inner stages with the local derivative operator i.e.
	\begin{equation}
	\begin{aligned}
		& \mathbf{u}_h^{(i)}= \mathbf{u}_h^n + \Delta t \sum_{j=1}^{i-1} a_{i j} \mathcal{G}(\mathbf{u}_{h}^{(j)} ) , \quad i = 1,2\hdots ,s,  \\
		\label{equation:cRKDG fully discrete}
		&\mathbf{u}_h^{n+1} = \mathbf{u}_h^n + \Delta t \sum_{i = 1}^{s} b_\ell \mathscr{F}(\mathbf{u}_{h}^{(i)} ).  \\  
	\end{aligned}
\end{equation}
We will say that a fully-discrete DG scheme is entropy bounded at any time level $t = t^n,$ if the following entropy principle is satisfied by it's solution 
\begin{align}
 s(\mathbf{u}_h( \mathbf{x}, t^n)) \ge s_0 \quad \text{and} \quad s_0 = \min_x s_0(\mathbf{x}), \
\end{align}
where $s_0(\mathbf{x})$ denotes the entropy of the initial condition. Alternatively, we will call a DG scheme entropy bounded if it's solution at any time level lies in the set $G.$ First, we illustrate the general methodology to construct an entropy bounded DG scheme for one-dimensional Euler equations. In order to do this, we consider the case when $\Omega$ is a bounded interval in $\mathbb{R}.$ We partition it into uniform subintervals each having length $h$. i.e.
\begin{align}
	\Omega = \bigcup_{i = 1}^K C_i, \quad C_i = [x_{i - \frac{1}{2}}, x_{i+\frac{1}{2}}] \quad \text{and} \quad x_{i+\frac{1}{2}} - x_{i-\frac{1}{2}} = h.
\end{align}
On each element $C_i $ we choose $N$-point Gauss-Lobatto quadrature so that it can exactly integrate the polynomials of degree at least $k.$ In other words, the value of $N$ should be large enough to satisfy the condition $ 2N - 3 \ge k.$ Denote these quadrature points on $C_i$ by 
\begin{align}
	\hat{S}_i := \{x_{i-\frac{1}{2}} = \hat{x}_{i}^{1}, \hat{x}_{i}^{2}, \ldots ,\hat{x}_{i}^{N} = x_{i+\frac{1}{2}}\},
\end{align}
The corresponding weights on the reference interval $[-\frac{1}{2}, \frac{1}{2}] $ will be denoted by $ \hat{w}_\alpha$ such that $\sum_{\alpha = 1}^{N}\hat{w}_\alpha = 1.$ In addition, we will use the notation $\mathbf{p}_i(x)$ to denote the polynomial solution of the DG method on the interval $C_i.$ Note that for entropy of the solution to be well defined at the quadrature points, the solution should preserve positivity of density and pressure at these points. We will say that $ \mathbf{u}_h \in G_0 $ on $\hat{S}_i$, if  $  \mathbf{u}_h(\hat{x}_i^{\alpha}) \in G_0$ for all $\hat{x}_i^{\alpha} \in\hat{S}_i.$ Similarly we will say that $ \mathbf{u}_h \in G_0 $ on $\tilde{S}_i$, if  $  \mathbf{u}_h(\hat{x}_i^{\alpha}) \in G_0$ for all $\hat{x}_i^{\alpha} \in\hat{S}_i$ and $\mathbf{u}^{+}_{i+\frac{1}{2}}, \mathbf{u}^{-}_{i-\frac{1}{2}} \in G_0.$
Let us now consider the DG scheme \eqref{equation:cRKDG fully discrete} with the following Butcher table .i.e the forward-Euler method 
\begin{align*}
\begin{array}{c|c}
	0 & 0\\
	\hline 
	 & 1\\
\end{array}
\end{align*}
Then the element averaged DG solution on any interval $C_i$ satisfies 
\begin{align}
	\bar{\mathbf{u}}_i^{n+1}= \bar{\mathbf{u}}_i^n - \lambda \left( \widehat{\mathbf{f}}(\mathbf{u}^{-}_{i+\frac{1}{2}}, \mathbf{u}^{+}_{i+\frac{1}{2}}) - \widehat{\mathbf{f}}(\mathbf{u}^{-}_{i-\frac{1}{2}}, \mathbf{u}^{+}_{i-\frac{1}{2}})\right),
\end{align}
where we omit the subscript $h$ for simplicity.
\begin{lemma}\textup{\cite{Zhang entropy, Lv Yu}} If $\mathbf{u}_h \in G_{0}$ on $\tilde{S}_i$ then the forward-Euler update $ \omega = \mathbf{u}_h + \Delta t \mathscr{F}(\mathbf{u}_h) $ satisfies the following entropy principle 
	\label{lemma:entropy rkdg}
	\begin{align}
		s(\bar{\omega} _i) \ge \min \{ \min _{1\le \alpha \le N} s( \mathbf{p}_i(\hat{x}_i^{\alpha})), s(\mathbf{u}^{+}_{i + \frac{1}{2}}), s(\mathbf{u}^{-}_{i - \frac{1}{2}})   \} ,
	\end{align}
	under the CFL constraint 
	\begin{align}
		\lambda \alpha_0 \le c_{\scalebox{0.7}{$ \mathscr{F}$} }.
		\label{equation: Cfl rkdg}
	\end{align}  
	\begin{proof}
		The proof of this result is easy therefore we omit it. Interested readers are encouraged to refer \cite{ Zhang entropy, Lv Yu} for the details of the proof.
	\end{proof}
\end{lemma}
The preceding Lemma implies that under the time step constraint \eqref{equation: Cfl rkdg} the following result holds
\begin{align}
	s(\bar{\mathbf{u}} _i^{n+1}) \ge \min \{ \min _{1\le \alpha \le N} s( \mathbf{p}_i(\hat{x}_i^{\alpha})), s(\mathbf{u}^{n,+}_{i + \frac{1}{2}}), s(\mathbf{u}^{n, -}_{i - \frac{1}{2}})   \}.
	\label{intermediate: entropy bounded}
\end{align}
If the solution at time level $t = t^n$ is entropy bounded and $s(\mathbf{u}^n_i) \ge s_0 $ on $ \tilde{S}_i$ then \eqref{intermediate: entropy bounded} implies that $s(\bar{\mathbf{u}} _i^{n+1}) \ge s_0.$ Now, one can apply the following entropy bound limiter to enforce 
$s(\mathbf{u}^{n+1}_i) \ge s_0 $ on $ \hat{S}_i$
\begin{align}
\mathbf{u}_i^{\mathcal{L}_{EB}} = \mathbf{u}_i + \theta_i^{EB} (\bar{\mathbf{u}}_i - \mathbf{u}_i) \quad \text{and} \quad \theta_i^{EB} = \frac{\tau}{\tau - [p(\bar{\mathbf{u}}_i) - \exp(s_0) \rho ^{\gamma}(\bar{\mathbf{u}}_i)]},
\label{limiter:entropy bound}
\end{align}
where $\tau = \min \{0, \min_{1\le \alpha \le N} \{p(\mathbf{u}_i(\hat{x}^{\alpha}_{i})) - \exp(s_0) \rho ^{\gamma}(\mathbf{u}_i(\hat{x}^{\alpha}_{i}))\}\}.$ This limiter does not destroy the accuracy and conservation of the original scheme. We refer the reader to \cite{Lv Yu} for more details about this limiter. For a higher order RK method in time the fully discrete DG scheme in \eqref{equation:RKDG fully discrete} is generally expressed in the Shu-Osher form as follows
\begin{align}
& \mathbf{u}_h^{(0)} = \mathbf{u}_h^{n},  \nonumber \\
& \mathbf{u}_h^{(i)} = \sum_{j = 0}^{i-1}\left(\alpha_{ij} \mathbf{u}_h^{(j)} + \Delta t \beta_{ij} \mathscr{F}(\mathbf{u}_h^{(j)})\right), \quad \quad 1\le i \le s
\label{equation:rkdg shu osher} \\
&  \mathbf{u}_h^{n+1} =  \mathbf{u}_h^{(s)}  \nonumber,
\end{align}
where at every stage $\sum_{j = 0}^{i-1} \alpha_{ij} = 1$ is satisfied. If we assume that the solution at $n^{\text{th}}$ time level or $ \mathbf{u}_h^{n}$ lies in $G$ and the positivity preserving limiter is applied after each stage of RK \eqref{equation:rkdg shu osher}, then a subsequent application of entropy bound limiter \eqref{limiter:entropy bound} at each stage implies that $\mathbf{u}_h^{n+1} \in G.$
\section{Entropy bounded cRKDG method}
In this section, we first introduce some preliminary results which will be useful for constructing entropy bounded cRKDG method for Euler equations later on. Before beginning, we recall some basic results associated with the Lax-Friedrichs scheme. Let us first define the Lax-Friedrichs average of two admissible states $U_{\ell}, U_r $ in $G_0$ with parameter $\alpha$ as follows:
\begin{align}
H(U_{\ell}, U_r; \alpha) := \frac{U_{\ell} + U_r}{2} - \frac{f(U_r) - f(U_l)}{2\alpha}
\label{equation: Riemann average}
\end{align}
Next, we consider a Riemann problem, $RP (U_{\ell}, U_r)$ for \eqref{equation:conservation law} with the following initial condition:
\begin{equation}
\mathbf{u}_0(x) = 
\begin{cases}
U_{\ell} & \, x<0, \\
U_r    & \, x > 0,
\end{cases}
\end{equation}
If a self-similar entropy solution to $RP (U_{\ell}, U_r)$ exists then one can establish a two point entropy inequality for the average in \eqref{equation: Riemann average} using the Reynolds transport theorem and repeated application of self-similarity.
\begin{lemma}\textup{(Kailiang Wu, \cite{kailiang wu})} Let $ U_{\ell} \, \text{and} \, U_r $ be two admissible states in $ G_0.$  Assume that the problem $RP(U_{\ell}, U_r) $ admits a self-similar entropy solution whose values lie in $G_0$ and 
	\begin{align*}
		\mathbf{u}(x, t) = \mathcal{R}(\xi) \quad \text{where} \quad \xi = x/t. 
	\end{align*} 
In addition, assume that all waves lie entirely in the region $|\xi | \le 
\alpha ,$ then we have
\begin{align}
\mathcal{U}\left(H(U_{\ell}, U_r; \alpha)\right) \leq \frac{\mathcal{U}(U_\ell)+ \mathcal{U}(U_r)}{2} - \frac{\mathcal{F}(U_r) - \mathcal{F}(U_\ell)}{2\alpha}.
\end{align}
Moreover, 
\begin{align*}
 H(U_{\ell}, U_r; \alpha) = \frac{1}{2\alpha} \int_{-\alpha}^{\alpha} \mathcal{R}(\xi) \, d \xi \quad \text{and} \quad H(U_{\ell}, U_r; \alpha) \in G_0.
\end{align*}
\label{lemma:entropy stability}
\end{lemma}
This result was utilized by the authors to establish weak entropy stability of the DG scheme with forward-Euler time stepping. We will use this result for proving entropy boundedness of the Lax-Friedrichs scheme when applied to the Euler equations. We highlight that it has already been established by Tadmor in \cite{Tadmor} using the entropy inequality proven by Lax \cite{Lax}. However, the Lax's proof of entropy inequality is based on homotopy method where the time step constraint depends upon the ratio of minimum and maximum eigenvalues of the Hessian matrix $\mathcal{U}^{''}$and therefore impractical for computational purpose. To ensure entropy boundedness of the LxF scheme under a suitable time step restriction we consider the following scheme for the Euler equations in 1d:
\begin{align}
U^{n+1}_i = \frac{U^{n}_{i-1} + U^{n}_{i+1}}{2} - \lambda \frac{\mathbf{f}(U_{i+1}^{n}) - \mathbf{f}(U_{i-1}^{n})}{2} \qquad \quad \lambda = \frac{\Delta t}{h},
\label{equation:Lax_friedrichs scheme}
\end{align}
where, $U_{i}^n \approx \frac{1}{h}\int _{C_i} \mathbf{u}(x, t^n) \, dx$ denotes the approximation of cell average at $t = t^n.$  
\begin{lemma}
	If the hypotheses of Lemma \ref{lemma:entropy stability} holds then the Lax-Friedrichs scheme \eqref{equation:Lax_friedrichs scheme}, for the Euler equations in 1d is positivity preserving under the CFL constraint
	\begin{align}
\lambda \alpha_0 \le 1.
\label{lemma2part:time constraint}
	\end{align}
In other words, if for every $i$, $U_i^{n} \in G_0 $ then $U_i^{n+1}\in G_0 $, provided \eqref{lemma2part:time constraint} holds. In addition, it satisfies the following entropy principle:
\begin{align}
s(U^{n+1}_{i}) \ge \min \{s(U^n_{i-1}), s(U^n_{i+1})\} \quad \quad \forall \, i.
\label{lemma2part:entropy principle}
\end{align}
\label{lemma2: entropy bounded Lax Freiedrichs}
\begin{proof}
Although, Lemma \ref{lemma:entropy stability} directly implies the positivity of the scheme, we still provide an alternative proof. Note that the scheme \eqref{equation:Lax_friedrichs scheme} can be rewritten as follows:
\begin{align*}
 U_i^{n+1} = \frac{1}{2} (1 - \lambda \alpha_0)\left( U_{i-1}^n + U_{i+1}^n\right) + \frac{\lambda \alpha_0}{2} \left( U_{i-1}^n - \alpha_0 ^{-1}\mathbf{f}( U_{i-1}^n)\right) + \frac{\lambda \alpha_0}{2} \left( U_{i+1}^n + \alpha_0 ^{-1}\mathbf{f}( U_{i+1}^n)\right).
\end{align*}
By the Lax-Friedrichs splitting property\cite{bound preserving crkdg} $ w \pm \alpha_0 ^{-1} \mathbf{f}(w)$ belong to $G_0,$ whenever $w$ is in $G_0$. Therefore, $ U_i^{n+1}$ being a convex combination of admissible states also belongs to $G_0$. Now the density and pressure of $ U_i^{n+1}$ are positive and using the similar arguments as in \cite{Tadmor}, one can easily get the entropy principle \eqref{lemma2part:entropy principle}.
\end{proof}
\end{lemma}
In the next result, we will use the entropy boundedness of the LxF scheme and  properties of Gauss-Lobatto quadrature.
\begin{lemma}
If $\mathbf{u}_h \in G_{0}$ on $\hat{S}_i$ then the forward-Euler updates $ \omega = \mathbf{u}_h \pm \Delta t \mathcal{G}(\mathbf{u}_h) $ satisfies the following entropy principle 
\label{lemma:entropy principle for crkdg}
\begin{align}
 s(\bar{\omega} _i) \ge  \min _{1\le \alpha \le N}\{ s( \mathbf{p}_i(\hat{x}_i^{\alpha})) \}, 
 \label{inequality: entropy inequality crkdg}
\end{align}
under the CFL constraint 
\begin{align}
\lambda \alpha_0 \le  c_{\scalebox{0.7}{$ \mathcal{G}$} .}
\label{equation: Cfl crkdg}
\end{align}  
\begin{proof}
From the weak bound preserving property in Lemma 3.3 of \cite{bound preserving crkdg}, $\bar{\omega}_{i} \in G_0.$ Therefore, the density and pressure of $\bar{\omega}_{i}$ are positive and $s(\bar{w} _i)$ is well defined. For simplicity we omit the subscript $i \, \text{and} \, h$ whenever possible. We will prove this Lemma for the case when $ \omega = \mathbf{u}_h + \Delta t \mathcal{G}(\mathbf{u}_h).$  Similar arguments can be applied to conclude the result for $ \omega = \mathbf{u}_h - \Delta t \mathcal{G}(\mathbf{u}_h).$
\begin{align}
\bar{\omega}  = \bar{\mathbf{u}} -\lambda (\mathbf{f}(\mathbf{u}_{i+\frac{1}{2}}^{-}) - \mathbf{f}(\mathbf{u}_{i-\frac{1}{2}}^{+})),
\end{align}
where $ \mathbf{u}_{i+\frac{1}{2}}^{-} = \mathbf{u}_h(x_{i+\frac{1}{2}}^{-})$ and $ \mathbf{u}_{i-\frac{1}{2}}^{-} = \mathbf{u}_h(x_{i-\frac{1}{2}}^{+}).$ The exactness of quadrature rule implies,
\begin{align}
\bar{\omega} = \sum_{q = 1}^{N} \hat{w}_q \mathbf{u}_{i, q} -\lambda (\mathbf{f}(\mathbf{u}_{i+\frac{1}{2}}^{-}) - \mathbf{f}(\mathbf{u}_{i-\frac{1}{2}}^{+})),
\end{align}
where we denote $ \mathbf{u}_{i, q} = \mathbf{u}(\hat{x}_{i}^q).$ Note that the choice of quadrature rule implies the following relation
\begin{align}
\bar{\omega} = \sum_{q = 2}^{N-1} \hat{w}_q \mathbf{u}_{i, q} + \hat{w}_1 H_1 + \hat{w}_N H_1,
\end{align}
where 
\begin{align}
 H_1 = \frac{\mathbf{u}_{i, 1} + \mathbf{u}_{i, N}}{2} - \frac{\lambda}{2\hat{w}_1}(\mathbf{f}(\mathbf{u}_{i, 1} ) - \mathbf{f}(\mathbf{u}_{i, N} )) .
\end{align}
Using Lemma \ref{lemma2: entropy bounded Lax Freiedrichs} we get the following inequality under the CFL constraint \eqref{equation: Cfl crkdg}  
\begin{align}
s(H_1) \ge \min \{ s(\mathbf{u}_{i, 1}), s(\mathbf{u}_{i, N})\}.
\end{align}
Next we use the quasi-concavity of the specific entropy to get 
\begin{align*}
s(\bar{\omega}) \ge \min \{\min _{2\le \alpha \le N-1} s( \mathbf{p}_i(\hat{x}_i^{\alpha})) , s(H_1)\},
\end{align*}
which in turn implies the result \eqref{inequality: entropy inequality crkdg}.
\end{proof}
\end{lemma}
\begin{lemma}
Assume that the $i^{\text{th}}$ stage of the cRKDG method in \eqref{equation:cRKDG fully discrete} can be expressed as the convex combination of the forward Euler steps $ \mathbf{u}_h^{(j)} + \Delta t \mu_{ij}^{\mathscr{F}} \mathscr{F}(\mathbf{u}_h^{(j)}) $ and $ \mathbf{u}_h^{(j)} \pm \Delta t \mu_{ij}^{\pm \mathcal{G}} \mathcal{G}(\mathbf{u}_h^{(j)}) $ where $j < i$ and $ \mu_{ij}^{\mathscr{F}}, \mu_{ij}^{\pm \mathcal{G}} > 0$. In addition, if we assume that for each $j < i$ we have $ \mathbf{u}^{(j)}_h \in G$ on $ \tilde{S}_K,$ then under the time step constraint 
\label{lemma: enrtopy boundedness for ith stage}
\begin{align}
\lambda \alpha_0 \le \min_{j < i} \left(\frac{c_{\scalebox{0.7}{$ \mathscr{F}$} }}{\mu_{ij}^{\mathscr{F}}}, \frac{c_{\scalebox{0.7}{$ \mathcal{G}$} }}{\mu_{ij}^{\pm \mathcal{G}}}\right)
\end{align}
$ \bar{\mathbf{u}}^{(i)}_h \in G $ on $ \hat{S}_K .$
\begin{proof}
We will omit the symbols $K$ and $h$ whenever it is possible. Since $ \mathbf{u}^{(j)} \in G $ on $ \tilde{S}_K $ therefore each $ \mathbf{u}^{(j)} \in G_0$ on $ \tilde{S}_K $ and by the hypotheses of this Lemma $ \bar{\mathbf{u}}^{(i)} \in G_0 $ on $ \hat{S}_K $ (see Lemma 3.4 in \cite{bound preserving crkdg}.) So, we just need to show that $ s(\bar{\mathbf{u}}^{(i)}) \ge s_0$ on $ \hat{S}_K $. We can write $ \bar{\mathbf{u}}^{(i)}$ as follows:
\begin{align}
\bar{\mathbf{u}}^{(i)} = \sum_{j = 1}^{i-1} \beta_{ij}^{\mathscr{F}}
\left( \overline{\mathbf{u}^{(j)} + \Delta t \mu_{ij}^{\mathscr{F}} \mathscr{F}(\mathbf{u}^{(j)})}\right) + \sum_{j = 1}^{i-1}\beta_{ij}^{\pm \mathcal{G}}  \left(\overline{\mathbf{u}^{(j)} \pm \Delta t \mu_{ij}^{\pm \mathcal{G}} \mathcal{G}(\mathbf{u}^{(j)})}\right), 
\end{align}
where $ \sum_{j = 1}^{i - 1} \beta_{ij}^{\mathscr{F}} + \beta_{ij}^{\pm \mathcal{G}}  = 1 \, \text{and} \, \beta_{ij}^{\mathscr{F}},\beta_{ij}^{\pm \mathcal{G}} \ge 0  .$ By using quasi-concavity of entropy we infer that
\begin{align}
s(\bar{\mathbf{u}}^{(i)}) \ge \min_{j < i} \{ \min \{ s( \overline{\mathbf{u}^{(j)} + \Delta t \mu_{ij}^{\mathscr{F}} \mathscr{F}(\mathbf{u}^{(j)})_K}), s( \overline{\mathbf{u}^{(j)} \pm \Delta t \mu_{ij}^{\pm \mathcal{G}} \mathcal{G}(\mathbf{u}^{(j)})_K})\}\}.
\label{lemmapart4: entropy of cell average}
\end{align}
Then a simple application of Lemma \ref{lemma:entropy rkdg} and Lemma \ref{lemma:entropy principle for crkdg} yields the following inequalities
\begin{align}
s( \overline{\mathbf{u}^{(j)} + \Delta t \mu_{ij}^{\mathscr{F}} \mathscr{F}(\mathbf{u}^{(j)})_K}) \ge \min \{ \min _{1\le \alpha \le N} s( \mathbf{p}^{(j)}_K(\hat{x}_K^{\alpha})), s(\mathbf{u}^{(j),+}_{K + \frac{1}{2}}), s(\mathbf{u}^{(j),-}_{K - \frac{1}{2}})   \},  
\label{lemmapart4: subsidiary equation 1}\\
s( \overline{\mathbf{u}^{(j)} \pm \Delta t \mu_{ij}^{\pm \mathcal{G}} \mathcal{G}(\mathbf{u}^{(j)})_K}) \ge  \min \{ \min _{1\le \alpha \le N} s( \mathbf{p}^{(j)}_K(\hat{x}_K^{\alpha})), s(\mathbf{u}^{(j),+}_{K + \frac{1}{2}}), s(\mathbf{u}^{(j),-}_{K - \frac{1}{2}})   \},
\label{lemmapart4: subsidiary equation 2}
\end{align}
where we denote the polynomial solution of $ \mathbf{u}^{(j)}$ on element $K$ 
by $ {\mathbf{p}}^{(j)}_K.$ Note that we also have $ s(\mathbf{u}^{(j)}) \ge s_0  $ on $ \tilde{S}_K$ which together with \eqref{lemmapart4: subsidiary equation 1}, \eqref{lemmapart4: subsidiary equation 2} and \eqref{lemmapart4: entropy of cell average} implies 
\begin{align}
s(\bar{\mathbf{u}}^{(i)}) \ge s_0 \qquad \text{on} \ &  \hat{S}_K .
\end{align}
Hence the conclusion follows.
\end{proof}
\end{lemma}	
\begin{theorem}
If each stage of the cRKDG method in \eqref{equation:cRKDG fully discrete} can be written as convex combination of forward-Euler steps of previous stages values with the operators $\pm \mathcal{G}, \mathscr{F}$, and bound preserving limiter along with entropy bound limiter is applied after each stage of RK then under the following time step constraint 
\begin{align}
\lambda \alpha_0 \le \min_{i, j} \left(\frac{c_{\scalebox{0.7}{$ \mathscr{F}$} }}{\mu_{ij}^{\mathscr{F}}}, \frac{c_{\scalebox{0.7}{$ \mathcal{G}$} }}{\mu_{ij}^{\pm \mathcal{G}}}\right),
\end{align}
$\mathbf{u}_h^{n} \in G $ on $ \hat{S}_K $ for all $K$ implies that $\mathbf{u}_h^{n+1}  \in G$ on $ \hat{S}_K  $ for all $K.$
\label{theorem:crkdg entropy boundedness}
\begin{proof}
The hypotheses of this Theorem already implies that $\mathbf{u}_h^{n+1}  \in G_0$ on $ \hat{S}_K  $ for all $K$ (see \cite{bound preserving crkdg}). So, we just need to show that $s(\mathbf{u}^{n+1}) \ge s_0$ on $ \hat{S}_K $ for every $K.$ Note that the cRKDG scheme with bound preserving limiter and subsequent application of entropy bound limiter after each RK stage can be written as follows:
	\begin{equation}
	\begin{aligned}
		& \mathbf{u}_h^{(i)}= \mathbf{u}_h^{n, \star} + \Delta t \sum_{j=1}^{i-1} a_{i j} \mathcal{G}(\mathbf{u}_{h}^{(j), \star} ) , \quad  \mathbf{u}_h^{(i), \star} = \mathcal{L}_{EB} \circ \mathcal{L}_{BP}(\mathbf{u}_h^{(i))}) \\
		\label{equation:cRKDG fully discrete with limiters}
		&\mathbf{u}_h^{n+1} = \mathbf{u}_h^{n, \star} + \Delta t \sum_{i = 1}^{s} b_\ell \mathscr{F}(\mathbf{u}_{h}^{(i), \star} ), \quad \mathbf{u}_h^{n+1, \star} = \mathcal{L}_{EB} \circ \mathcal{L}_{BP}(\mathbf{u}_h^{n+1}) \\  
	\end{aligned}
\end{equation}
From now on wards we will omit the symbol $h$ for convenience. For the inner stages of the cRKDG method with limiters, each of the stage values can be written as 
\begin{align}
	\mathbf{u}^{(i)} = \sum_{j=1}^{i-1} \beta_{ij} \left(\mathbf{u}^{(j)} \pm \mu_{ij}^{\pm \mathcal{G}} \mathcal{G}(\mathbf{u}^{(j), \star})\right) \quad \text{and} \quad \sum_{j = 1}^{i-1} \beta_{ij} = 1.
\end{align}
The subsequent application of limiters does not alter the cell average of the solution therefore we have 
\begin{align*}
	\bar{\mathbf{u}}^{(i)} = \sum_{j=1}^{i-1} \beta_{ij} \left(\bar{\mathbf{u}}^{(j)} \pm \mu_{ij}^{\pm \mathcal{G}} \overline{\mathcal{G}(\mathbf{u}^{(j), \star})_K}\right) = \sum_{j=1}^{i-1} \beta_{ij} \overline{\left(\mathbf{u}^{(j), \star} \pm \mu_{ij}^{\pm \mathcal{G}} \mathcal{G}(\mathbf{u}^{(j), \star})_K\right)} 
\end{align*}
Then by Lemma \ref{lemma: enrtopy boundedness for ith stage}, $ 	\bar{\mathbf{u}}^{(i)} \in G$ on $ \hat{S}_K, $ which in turn implies $\mathbf{u}^{(i), \star} \in G$ on $ \hat{S}_K.$ Similarly, the final stage of the method will be of the form 
\begin{align*}
\mathbf{u}^{n+1} = \sum_{i = 1}^{s} \beta_{i}^{\mathscr{F}}\left( \mathbf{u}^{(i)} + \Delta t \mu_{i}^{\mathscr{F}} \mathscr{F}(\mathbf{u}^{(i), \star})\right) + \sum_{i = 1}^{s} \beta_{i}^{\pm \mathcal{G}} \left(\mathbf{u}^{(i)} \pm \Delta t \mu_{i}^{\pm \mathcal{G}} \mathcal{G}(\mathbf{u}^{(i), \star})\right),
\end{align*} 
where $ \sum_{i = 1}^{s} \beta_{i}^{\mathscr{F}} +  \beta_{i}^{\mathcal{\pm G}} = 1 $ and $ \beta_{i}^{\mathscr{F}}, \beta_{i}^{\mathcal{\pm G}} \ge 0 .$ Therefore, the cell average on element $K$ is
\begin{align*}
\bar{\mathbf{u}}^{n+1} & = \sum_{i = 1}^{s} \beta_{i}^{\mathscr{F}} \overline{\left( \mathbf{u}^{(i)} + \Delta t \mu_{i}^{\mathscr{F}} \mathscr{F}(\mathbf{u}^{(i), \star})_K\right)} + \sum_{i = 1}^{s} \beta_{i}^{\pm \mathcal{G}} \overline{\left(\mathbf{u}^{(i)} \pm \Delta t \mu_{i}^{\pm \mathcal{G}} \mathcal{G}(\mathbf{u}^{(i), \star})_K\right)} \\
  & = \sum_{i = 1}^{s} \beta_{i}^{\mathscr{F}} \overline{\left( \mathbf{u}^{(i),\star} + \Delta t \mu_{i}^{\mathscr{F}} \mathscr{F}(\mathbf{u}^{(i), \star})_K\right)} + \sum_{i = 1}^{s} \beta_{i}^{\pm \mathcal{G}} \overline{\left(\mathbf{u}^{(i), \star} \pm \Delta t \mu_{i}^{\pm \mathcal{G}} \mathcal{G}(\mathbf{u}^{(i), \star})_K\right)}.
\end{align*}
Which yields that $ \mathbf{u}^{n+1} \in G $ on $\hat{S}_K$ for all $K.$
\end{proof}
\end{theorem}
\section{Two dimensional case}
In this section we consider the two-dimensional Euler equations given by
\begin{align}
	\mathbf{u}_t + \nabla \cdot \mathbf{F}(\mathbf{u}) &= 0
	\quad \text{and} \quad
	\mathbf{F}(\mathbf{u})
	= \left(\mathbf{f}_1(\mathbf{u}),\mathbf{f}_2(\mathbf{u})\right),
	\label{eq:euler}
\end{align}
and
\begin{align}
	\mathbf{u} =
	\begin{pmatrix}
		\rho \\
		\rho v_1 \\
		\rho v_2 \\
		\rho e
	\end{pmatrix},
\quad
	\mathbf{f}_1(\mathbf{u}) =
	\begin{pmatrix}
		\rho v_1 \\
		\rho v_1^2 + p \\
		\rho v_1v_2 \\
		(\rho e+p)v_1
	\end{pmatrix},
\quad
	\mathbf{f}_2(\mathbf{u}) =
	\begin{pmatrix}
		\rho v_2 \\
		\rho v_1v_2 \\
		\rho v_2^2 + p \\
		(\rho e+p)v_2
	\end{pmatrix}.
	\label{eq:euler_variables_fluxes}
\end{align}
where $v_1$ and $v_2$ are components of velocity $ \mathbf{v}$ in $x$ and $y$ directions respectively. \\
We assume that we have a uniform rectangular mesh and we denote the cell $(i, j)$ by 
\begin{align}
	I_{ij} = [x_{i-\frac{1}{2}}, x_{i+\frac{1}{2}}] \times [y_{j-\frac{1}{2}}, y_{i+\frac{1}{2}}],\, \Delta x = x_{i+\frac{1}{2}} - x_{i-\frac{1}{2}}, \, 
	\Delta y = y_{j+\frac{1}{2}} - y_{j-\frac{1}{2}}.
	\nonumber
\end{align} 
Let $\mathbf{p}^n_{ij}(x, y)$ be the vector of polynomial approximation of the solution of DG scheme at $n^ {\text{th}}$ time level. Then the cell average of DG solution on element $I_{ij}$ would be 
\begin{align}
 \bar{\mathbf{u}}_{ij}^n = \frac{1}{\Delta x \Delta y}\int_{I_{ij}} \mathbf{p}^n_{ij}(x,\, y) \, dx \, dy.
\end{align}
Let $ \{x_i^{q}\}_{q = 1}^{L} $, $ \{y_j^{q}\}_{q = 1}^{L}   $  be an $L$-point Gauss-Quadrature and $\{ w_q\}_{q = 1}^L$ be the corresponding weights on the reference interval $[ -\frac{1}{2}, \frac{1}{2}].$ We also choose a $N-$point Gauss-Lobatto quadrature $ \{\hat{x}_i^{\alpha}\}_{\alpha = 1}^{N} $, $ \{\hat{y}_j^{\alpha}\}_{\alpha = 1}^{N}$ and weights $\{\hat{w}_\alpha\}_{\alpha=1}^{N}.$ On every element $I_{ij}$ we will use the following set of quadrature points
\begin{align}
S_{ij} := (\hat{S}_i^x \times \bar{S}_j^y) \cup (\bar{S}_i^x \times \hat{S}_j^y), 
\label{quadrature}
\end{align}
where $ \hat{S}_i^x := \{\hat{x}_i^{\alpha}\}_{\alpha = 1}^{N}, \bar{S}_j^y := \{y_j^{q}\}_{q = 1}^{L}, \bar{S}_i^x := \{x_i^{q}\}_{q = 1}^{L}, \hat{S}_j^y := \{\hat{y}_j^{\alpha}\}_{\alpha = 1}^{N}.$
The traces of the solution on the edges of the cell $(i, j)$ will be denoted by 
\begin{align*}
	\mathbf{u}^{\mp}_{i\pm\frac{1}{2}, j}(y) = \mathbf{p}_{ij}(x_{i\pm\frac{1}{2}},\, y), \,	\mathbf{u}^{\mp}_{i, j\pm\frac{1}{2}}(x) = \mathbf{p}_{ij}(x,\, y_{j\pm\frac{1}{2}}).
\end{align*}
Similar to the case of one-dimension, we will say that $\mathbf{u}_h \in G$ on $S_{ij}$ if $ \mathbf{p}_{ij}(\hat{x}_i^{\alpha},  y_j^{q}), \mathbf{p}_{ij}(x_i^{q}, \hat{y}_j^{\alpha}) \in G $ for all $\alpha = 1, 2,.\ldots, N$ and $q = 1, 2,.\ldots, L.$
In addition, $\mathbf{u}_h \in G$ on $\tilde{S}_{ij}$ if $\mathbf{u}_h \in G$ on $S_{ij}$ and $ \mathbf{u}^{\pm}_{i \pm \frac{1}{2}, q},  \mathbf{u}^{\pm}_{q, j \pm \frac{1}{2}} \in G $ for all $q = 1, 2,.\ldots, L.$
 The DG scheme with forward-Euler update can be written as follows
\begin{align}
\bar{\mathbf{u}}_{ij}^{n+1}  = \bar{\mathbf{u}}_{ij}^{n} & - \frac{\Delta t}{\Delta x \Delta y } \int_{y_{j-\frac{1}{2}}}^{y_{j+ \frac{1}{2}}}\left(\hat{\mathbf{f}}_1(\mathbf{u}_{i+\frac{1}{2}, j}^{-}(y), \mathbf{u}_{i+\frac{1}{2}, j}^{+}(y)) - \hat{\mathbf{f}}_1(\mathbf{u}_{i-\frac{1}{2}, j}^{-}(y), \mathbf{u}_{i-\frac{1}{2}, j}^{+}(y))  \right) \, dy \nonumber \\
& - \frac{\Delta t}{\Delta x \Delta y} \int_{x_{i-\frac{1}{2}}}^{x_{i+ \frac{1}{2}}}\left(\hat{\mathbf{f}}_2(\mathbf{u}_{i, j+\frac{1}{2}}^{-}(x), \mathbf{u}_{i, j+\frac{1}{2}}^{+}(x)) - \hat{\mathbf{f}}_2(\mathbf{u}_{i, j-\frac{1}{2}}^{-}(x), \mathbf{u}_{i, j-\frac{1}{2}}^{+}(x))  \right) \, dx,
\label{mathequation: 2d DG scheme} 
\end{align}
where $\hat{\mathbf{f}}_1(.,.)$ and $\hat{\mathbf{f}}_2(.,.) $ are one-dimensional numerical fluxes. In the present scenario, we use the Lax-Friedrichs numerical flux defined by 
\begin{align*}
\hat{\mathbf{f}}_1(\mathbf{a}, \mathbf{b}) = \frac{1}{2} [\mathbf{f}_1(\mathbf{a}) + \mathbf{f}_1(\mathbf{b}) - \alpha_1 (\mathbf{b} - \mathbf{a})]  \qquad \alpha_1 = \| (| v_1| + c) \|_{\infty}, \\
\hat{\mathbf{f}}_2(\mathbf{a}, \mathbf{b}) = \frac{1}{2} [\mathbf{f}_2(\mathbf{a}) + \mathbf{f}_2(\mathbf{b}) - \alpha_2 (\mathbf{b} - \mathbf{a})]  \qquad \alpha_2 = \| (| v_2| + c) \|_{\infty}.
\end{align*}
The cell average of the solution can be decomposed by using the Jiang-Liu convex decomposition \cite{Jiang Liu} in equation \eqref{quadrature} to get
\begin{align}
&\bar{\mathbf{u}}_{ij}^n = \frac{1}{2} \sum_{\alpha = 1}^{N} \sum_{q = 1}^{L} \hat{w}_{\alpha} w_{q} \, \mathbf{p}^{n}_{ij}(\hat{x}_i^{\alpha}, y_j^{q}) + \frac{1}{2} \sum_{\alpha = 1}^{N} \sum_{q = 1}^{L} \hat{w}_{\alpha} w_{q} \, \mathbf{p}^{n}_{ij}(x_i^{q}, \hat{y}_j^{\alpha}), 
\label{eq:decomposed cell average}\\
\nonumber 
\end{align}
 There are several other choices like the the general convex decomposition or the optimal convex decomposition\cite{S cui, S cui kailiang wu}. 
 Using equation \eqref{eq:decomposed cell average}, the cell average update in \eqref{mathequation: 2d DG scheme} can be written as follows
\begin{align}
\bar{\mathbf{u}}_{ij}^{n+1}   & = \bar{\mathbf{u}}_{ij}^{n} - \lambda_1 \sum_{q = 1}^{L} w_q \left(\hat{\mathbf{f}}_1(\mathbf{u}^{-}_{i + \frac{1}{2}, q},\mathbf{u}^{+}_{i + \frac{1}{2}, q} ) - \hat{\mathbf{f}}_1(\mathbf{u}^{-}_{i - \frac{1}{2}, q},\mathbf{u}^{+}_{i - \frac{1}{2}, q} )\right) 
  \label{eq:2d dg decomposed}
 \nonumber \\
 & \qquad - \lambda_2 \sum_{q = 1}^{L} w_q \left(\hat{\mathbf{f}}_2(\mathbf{u}^{-}_{q, j + \frac{1}{2}},\mathbf{u}^{+}_{q, j + \frac{1}{2}} ) - \hat{\mathbf{f}}_2(\mathbf{u}^{-}_{q, j - \frac{1}{2}},\mathbf{u}^{+}_{q, j - \frac{1}{2}} )\right) \nonumber \\
 & =  \frac{1}{2} \sum_{\alpha = 1}^{N} \sum_{q = 1}^{L} \hat{w}_{\alpha} w_{q} \, \mathbf{p}^{n}_{ij}(\hat{x}_i^{\alpha}, y_j^{q}) + \frac{1}{2} \sum_{\alpha = 1}^{N} \sum_{q = 1}^{L} \hat{w}_{\alpha} w_{q} \, \mathbf{p}^{n}_{ij}(x_i^{q}, \hat{y}_j^{\alpha})  \\
 & \qquad - \lambda_1 \sum_{q = 1}^{L} w_q \left(\hat{\mathbf{f}}_1(\mathbf{u}^{-}_{i + \frac{1}{2}, q},\mathbf{u}^{+}_{i + \frac{1}{2}, q} ) - \hat{\mathbf{f}}_1(\mathbf{u}^{-}_{i - \frac{1}{2}, q},\mathbf{u}^{+}_{i - \frac{1}{2}, q} )\right) \nonumber \\
 & \qquad - \lambda_2 \sum_{q = 1}^{L} w_q \left(\hat{\mathbf{f}}_2(\mathbf{u}^{-}_{q, j + \frac{1}{2}},\mathbf{u}^{+}_{q, j + \frac{1}{2}} ) - \hat{\mathbf{f}}_2(\mathbf{u}^{-}_{q, j - \frac{1}{2}},\mathbf{u}^{+}_{q, j - \frac{1}{2}} )\right).
 \nonumber
\end{align}
where we denote $ \mathbf{u}^{\mp}_{i \pm \frac{1}{2}, q} =  \mathbf{u}^{\mp}_{i \pm \frac{1}{2}}(y^q_{j}),  \mathbf{u}^{\mp}_{q, j \pm \frac{1}{2}} =  \mathbf{u}^{\mp}_{j \pm \frac{1}{2}}(x^q_{i}) \, \text{and} \, \lambda_1 = \Delta t / \Delta x, \lambda_2 = \Delta t / \Delta y.$ 
We also define one-dimensional slice averages as follows:
\begin{align}
	\bar{\mathbf{u}}_{ij, q}^{x, n} = \sum_{\alpha = 1}^{N} \hat{w}^{\alpha} \mathbf{p}_{ij}^{n}(\hat{x}_{i}^{\alpha}, y_{j}^q), \quad \bar{\mathbf{u}}_{ij, q}^{y, n} = \sum_{\alpha = 1}^{N} \hat{w}^{\alpha} \mathbf{p}_{ij}^{n}(x_{i}^{q}, \hat{y}_{j}^{\alpha}).
\end{align}
This implies that the equation \eqref{eq:decomposed cell average} can be rewritten in terms of the slice averages $	\bar{\mathbf{u}}_{ij, q}^{x, n}$ and $\bar{\mathbf{u}}_{ij, q}^{y, n}$, so that we have
$$
\begin{aligned}
	\bar{\mathbf{u}}_{ij}^n = \frac{1}{2} \sum_{q = 1}^{L} w_{q} \, \bar{\mathbf{u}}_{ij, q}^{x, n} + \frac{1}{2} \sum_{q = 1}^{L} w_{q} \, \bar{\mathbf{u}}_{ij, q}^{y, n}.
\end{aligned}
$$
\begin{lemma}\textup{(Kailiang Wu, \cite{kailiang wu})} For each unit vector $ \textbf{n} $ let us define 
	\begin{align}
H(U_{\ell}, U_r; \alpha) := \frac{U_{\ell} + U_r}{2} - \frac{\mathbf{F}(U_r)\cdot \mathbf{n} - \mathbf{F}(U_l)\cdot \mathbf{n}}{2\alpha}.
	\end{align}
	If the one-dimensional Riemann problem 
	\begin{align}
		\mathbf{u}_t + \left(\mathbf{F}(\mathbf{u})\cdot \mathbf{n}\right)_{\xi} = 0
	\end{align}
with left and right states $U_\ell, U_r \in G_0,$ admits an entropy solution whose values lie in $G_0$ and all waves lie entirely in the region $ |\xi| \le \alpha.$ Then 
\begin{align}
\mathcal{U}\left(H(U_{\ell}, U_r; \alpha)\right) \leq \frac{\mathcal{U}(U_\ell)+ \mathcal{U}(U_r)}{2} - \frac{\mathcal{F}(U_r) \cdot \mathbf{n} - \mathcal{F}(U_\ell) \cdot \mathbf{n}}{2\alpha}
\end{align} 
and $ H(U_{\ell}, U_r; \alpha) \in G_0.$
\label{lemma2:kailiang wu}
\end{lemma}
Next, we define the directional LxF averages on each cell $ (i, j)$
\begin{align}
 \mathcal{A}_{ij, q}^{x} &:= \frac{\mathbf{u}_{i - \frac{1}{2}, q}^{+}+ \mathbf{u}_{i + \frac{1}{2}, q}^{-}}{2} - \lambda_1\frac{ \mathbf{f}_{1}(\mathbf{u}_{i + \frac{1}{2}, q}^{-}) - \mathbf{f}_{1}(\mathbf{u}_{i - \frac{1}{2}, q}^{+})}{2}, \\
\mathcal{A}_{ij, q}^{y} &:= \frac{\mathbf{u}_{q, j - \frac{1}{2}}^{+}+ \mathbf{u}_{q, j + \frac{1}{2}}^{-}}{2} - \lambda_2\frac{ \mathbf{f}_{2}(\mathbf{u}_{q, j + \frac{1}{2}}^{-}) - \mathbf{f}_{2}(\mathbf{u}_{q, j - \frac{1}{2}}^{+})}{2}.
\end{align}
Using the preceding Lemma yields the following result.
\begin{lemma}
	We assume that the hypotheses of Lemma \ref{lemma2:kailiang wu} holds and 
	$ \mathbf{u}_{i \pm \frac{1}{2}, q}^{\mp}, \mathbf{u}_{q, j \pm \frac{1}{2}}^{\mp} \in G_{0}.$ If the CFL constraints $ \lambda_1 \alpha_0 \le 1 $ and $\lambda_2 \alpha_0 \le 1$ are satisfied then $\mathcal{A}_{ij, q}^{x}, \mathcal{A}_{ij, q}^{y} \in G_0$ and the following entropy principles are satisfied by the directional LxF averages 
	\begin{align}
  & s(\mathcal{A}_{ij, q}^{x} ) \ge \min\{s(\mathbf{u}_{i - \frac{1}{2}, q}^{+}), s(\mathbf{u}_{i + \frac{1}{2}, q}^{-}) \}, \\
  & s(\mathcal{A}_{ij, q}^{y} ) \ge \min\{s(\mathbf{u}_{q, j - \frac{1}{2}}^{+}), s(\mathbf{u}_{q, j + \frac{1}{2}}^{-}) \}.
	\end{align}
\begin{proof}
	The proof  follows essentially the same idea as outlined in Lemma \ref{lemma2: entropy bounded Lax Freiedrichs}, therefore we omit it.
\end{proof}
\end{lemma}
Next, we establish an entropy principle similar to Lemma \ref{lemma:entropy principle for crkdg} for the two-dimensional case.
\begin{lemma}
If $\mathbf{u}_h \in G_0$ on $S_{ij}$ and the following CFL constraints is satisfied   
\begin{align}
\lambda \alpha_0 \, \le \, \frac{ c_{\scalebox{0.7}{$\mathcal{G}$}}}{2} \quad \text{where} \quad \lambda = \max \{\lambda_1, \lambda_2\}.
\label{lemma5: CFl condition}
\end{align}
then the forward-Euler update $\omega = \mathbf{u}_h \pm \Delta t \mathcal{G}(\mathbf{u}_h)  $ satisfy the following entropy principle 
\begin{align}
	s(\bar{\omega}_{ij}) \ge \min_{\alpha, q} \{ s(\mathbf{p}_{ij}(\hat{x}^{\alpha}_i, y_j^q)), s(\mathbf{p}_{ij}(x_i^q, \hat{y}^{\alpha}_j)) \},
	\label{equation:entropy principle G}
\end{align}
\begin{proof}
Under the CFL constraint \eqref{lemma5: CFl condition}, the element averaged solution, $\bar{\omega}_{ij} \in G_0$ using the weak bound preserving property for the operators $\pm\mathcal{G}$ in \cite{bound preserving crkdg}. Therefore, $s(\bar{\omega}_{ij})$ is well defined. This property will be proved only for the case $\omega = \mathbf{u}_h + \Delta t \mathcal{G}(\mathbf{u}_h) $. For the other case, one can use similar arguments. Using the quadrature and the convex decomposition as in \eqref{eq:2d dg decomposed} we have 
\begin{align}
\bar{\omega}_{ij} & = \left(\overline{\textbf{u}_h + \Delta t \mathcal{G}(\mathbf{u}_h)}\right)_{I_{ij}} \nonumber \\
& = \bar{\mathbf{u}_{ij}} -\lambda_1 \left(\sum_{q = 1}^{L} \mathbf{f}_1(\mathbf{u}_{i +\frac{1}{2}, q}^{-}) - \mathbf{f}_1(\mathbf{u}_{i -\frac{1}{2}, q}^{+})\right) -\lambda_2 \left(\sum_{q = 1}^{L} \mathbf{f}_2(\mathbf{u}_{q, j +\frac{1}{2}}^{-}) - \mathbf{f}_2(\mathbf{u}_{q, j -\frac{1}{2}}^{+})\right) \nonumber \\
& =  \frac{1}{2} \sum_{\alpha = 2}^{N-1} \sum_{q = 1}^{L} \hat{w}_{\alpha} w_{q} \, \mathbf{p}^{n}_{ij}(\hat{x}_i^{\alpha}, y_j^{q}) + \frac{1}{2} \sum_{\alpha = 2}^{N-1} \sum_{q = 1}^{L} \hat{w}_{\alpha} w_{q} \, \mathbf{p}^{n}_{ij}(x_i^{q}, \hat{y}_j^{\alpha})  \nonumber \\
& \qquad + \frac{1}{2} \hat{w}_1 \sum_{q = 1}^{L} w_q ( \mathbf{u}_{i + \frac{1}{2}, q}^{-}  + \mathbf{u}_{i - \frac{1}{2}, q}^{+}) -\lambda_1 \left(\sum_{q = 1}^{L} \mathbf{f}_1(\mathbf{u}_{i +\frac{1}{2}, q}^{-}) - \mathbf{f}_1(\mathbf{u}_{i -\frac{1}{2}, q}^{+})\right) \nonumber \\
& \qquad + \frac{1}{2} \hat{w}_1 \sum_{q = 1}^{L} w_q ( \mathbf{u}_{q, j + \frac{1}{2}}^{-}  + \mathbf{u}_{q, j - \frac{1}{2}}^{+}) -\lambda_2 \left(\sum_{q = 1}^{L} \mathbf{f}_2(\mathbf{u}_{q, j +\frac{1}{2}}^{-}) - \mathbf{f}_2(\mathbf{u}_{q, j -\frac{1}{2}}^{+})\right) \nonumber \\
& =  \frac{1}{2} \sum_{\alpha = 2}^{N-1} \sum_{q = 1}^{L} \hat{w}_{\alpha} w_{q} \, \mathbf{p}^{n}_{ij}(\hat{x}_i^{\alpha}, y_j^{q}) + \frac{1}{2} \sum_{\alpha = 2}^{N-1} \sum_{q = 1}^{L} \hat{w}_{\alpha} w_{q} \, \mathbf{p}^{n}_{ij}(x_i^{q}, \hat{y}_j^{\alpha})  \nonumber \\
& \qquad + \hat{w}_1 \sum_{q = 1}^{L} w_q \left[ \frac{( \mathbf{u}_{i + \frac{1}{2}, q}^{-}  + \mathbf{u}_{i - \frac{1}{2}, q}^{+})}{2} - \frac{2 \lambda_1}{2 \hat{w}_1} \left(\mathbf{f}_1(\mathbf{u}_{i +\frac{1}{2}, q}^{-}) - \mathbf{f}_1(\mathbf{u}_{i -\frac{1}{2}, q}^{+})\right)\right] \nonumber \\
& \qquad + \hat{w}_1 \sum_{q = 1}^{L} w_q \left[\frac{( \mathbf{u}_{q, j + \frac{1}{2}}^{-}  + \mathbf{u}_{q, j - \frac{1}{2}}^{+})}{2} -\frac{2\lambda_2}{2 \hat{w}_1} \left( \mathbf{f}_2(\mathbf{u}_{q, j +\frac{1}{2}}^{-}) - \mathbf{f}_2(\mathbf{u}_{q, j -\frac{1}{2}}^{+})\right)\right]. 
\label{equation: convex decomposition of cell average 2d}
\end{align}
Equation \eqref{equation: convex decomposition of cell average 2d} implies that $ \bar{\omega}_{ij} $ being a convex combination of admissible states in $G_0$  also belongs to $ G_0$ under the time step constraint \eqref{lemma5: CFl condition}. Finally, the quasi-concavity of $s$ implies the entropy-principle \eqref{equation:entropy principle G}.
\end{proof}
\end{lemma}
\begin{lemma}
If $\mathbf{u}_h \in G_0$ on $\tilde{S}_{ij}$ then the forward-Euler update $ \omega = \mathbf{u}_h + \Delta t \mathscr{F}(\mathbf{u}_h) $ satisfy the following entropy principle
\begin{align}
s(\bar{\omega}_{ij}) \ge\, \min_{q}\{ \min_{\alpha, q} \{ s(\mathbf{p}_{ij}(\hat{x}^{\alpha}_i, y_j^q)), s(\mathbf{p}_{ij}(x_i^q, \hat{y}^{\alpha}_j)) \}, s(\mathbf{u}^{\pm}_{i \pm \frac{1}{2}, q}), s(\mathbf{u}^{\pm}_{q, j \pm \frac{1}{2}})\}, 
\end{align}
under the CFL constraint
\begin{align}
\lambda \alpha_0 \, \le \, \frac{ c_{\scalebox{0.7}{$\mathcal{F}$}}}{2}.
\end{align}
\begin{proof}
The element averaged solution $\bar{\mathbf{u}}_{ij}$ can be decomposed using the Jiang-Liu convex decomposition as in \eqref{eq:decomposed cell average}
to get 
\begin{align}
 \bar{\omega}_ {ij}  =  \frac{1}{2} & \sum_{\alpha = 2}^{N-1} \sum_{q = 1}^{L} \hat{w}_{\alpha} w_{q} \, \mathbf{p}^{n}_{ij}(\hat{x}_i^{\alpha}, y_j^{q}) + \frac{1}{2} \sum_{\alpha = 2}^{N-1} \sum_{q = 1}^{L} \hat{w}_{\alpha} w_{q} \, \mathbf{p}^{n}_{ij}(x_i^{q}, \hat{y}_j^{\alpha}) \nonumber \\
 & + \frac{1}{2}\sum_{q = 1}^{L}w_q \hat{w}_1 \left[ 
  \mathbf{u}_{i+ \frac{1}{2}, q}^{-} - \frac{2 \lambda_1}{\hat{w}_1}\left(\hat{\mathbf{f}}_1(\mathbf{u}_{i+\frac{1}{2}, q}^{-}, \mathbf{u}_{i+\frac{1}{2}, q}^{+}) - \hat{\mathbf{f}}_1(\mathbf{u}_{i-\frac{1}{2}, q}^{+}, \mathbf{u}_{i+\frac{1}{2}, q}^{-})\right)\right]
  \nonumber \\
  & + \frac{1}{2}\sum_{q = 1}^{L}w_q \hat{w}_1 \left[ 
  \mathbf{u}_{i- \frac{1}{2}, q}^{+} - \frac{2 \lambda_1}{\hat{w}_1}\left(\hat{\mathbf{f}}_1(\mathbf{u}_{i-\frac{1}{2}, q}^{+}, \mathbf{u}_{i+\frac{1}{2}, q}^{-}) - \hat{\mathbf{f}}_1(\mathbf{u}_{i-\frac{1}{2}, q}^{-}, \mathbf{u}_{i-\frac{1}{2}, q}^{+})\right)\right]
  \nonumber \\
   & + \frac{1}{2}\sum_{q = 1}^{L}w_q \hat{w}_1 \left[ 
  \mathbf{u}_{q, j + \frac{1}{2}}^{-} - \frac{2 \lambda_2}{\hat{w}_1}\left(\hat{\mathbf{f}}_2(\mathbf{u}_{q, j+\frac{1}{2}}^{-}, \mathbf{u}_{q, j+\frac{1}{2}}^{+}) - \hat{\mathbf{f}}_2(\mathbf{u}_{q, j-\frac{1}{2}}^{+}, \mathbf{u}_{q, j+\frac{1}{2}}^{-})\right)\right]
  \nonumber \\
  & + \frac{1}{2}\sum_{q = 1}^{L}w_q \hat{w}_1 \left[ 
  \mathbf{u}_{q, j-\frac{1}{2}}^{+} - \frac{2 \lambda_2}{\hat{w}_1}\left(\hat{\mathbf{f}}_1(\mathbf{u}_{q, j -\frac{1}{2}}^{+}, \mathbf{u}_{q, j+\frac{1}{2}}^{-}) - \hat{\mathbf{f}}_1(\mathbf{u}_{q, j-\frac{1}{2}}^{-}, \mathbf{u}_{q, j-\frac{1}{2}}^{+})\right)\right]
  \nonumber
\end{align}
Therefore, the cell average $\bar{\omega} _{ij}$ can be written as a convex combination of finite volume updates and admissible states in $ G_0$ as in \eqref{eq:2d dg decomposed}. The quasi-concavity of $s$ then yields the desired inequality.
\end{proof}
\end{lemma}
In the forthcoming results we will denote any element $ I_{r, s} $ by $K$ to avoid any possible confusion with the indices $i, j$ which will be used to denote the stage values.
\begin{lemma}
	Assume that the $i^{\text{th}}$ stage of the cRKDG method in \eqref{equation:cRKDG fully discrete} can be expressed as the convex combination of the forward Euler steps $ \mathbf{u}_h^{(j)} + \Delta t \mu_{ij}^{\mathscr{F}} \mathscr{F}(\mathbf{u}_h^{(j)}) $ and $ \mathbf{u}_h^{(j)} \pm \Delta t \mu_{ij}^{\pm \mathcal{G}} \mathcal{G}(\mathbf{u}_h^{(j)}) $ where $j < i$ and $ \mu_{ij}^{\mathscr{F}}, \mu_{ij}^{\pm \mathcal{G}} > 0$. In addition, if we assume that for each $j < i$ we have $ \mathbf{u}^{(j)}_h \in G$ on $ \tilde{S}_K,$ then under the time step constraint 
	\label{lemma: enrtopy boundedness for ith stage 2d}
	\begin{align}
		\lambda \alpha_0 \le \min_{j < i} \left(\frac{c_{\scalebox{0.7}{$ \mathscr{F}$} }}{2\mu_{ij}^{\mathscr{F}}}, \frac{c_{\scalebox{0.7}{$ \mathcal{G}$} }}{2\mu_{ij}^{\pm \mathcal{G}}}\right)
	\end{align}
	$ \bar{\mathbf{u}}^{(i)}_h \in G $ on $ S_K .$
	\begin{proof}
 The proof is similar to the proof of Lemma \ref{lemma: enrtopy boundedness for ith stage} therefore we omit it.
	\end{proof}
\end{lemma}
\begin{theorem}
	If each stage of the cRKDG method in \eqref{equation:cRKDG fully discrete} can be written as convex combination of forward-Euler steps of previous stages values with the operators $\pm \mathcal{G}, \mathscr{F}$, and bound preserving limiter along with entropy bound limiter is applied after each stage of RK then under the following time step constraint 
	\begin{align}
		\lambda \alpha_0 \le \min_{i, j} \left(\frac{c_{\scalebox{0.7}{$ \mathscr{F}$} }}{2\mu_{ij}^{\mathscr{F}}}, \frac{c_{\scalebox{0.7}{$ \mathcal{G}$} }}{2\mu_{ij}^{\pm \mathcal{G}}}\right),
		\label{Cfl: crkdg 2d}
	\end{align}
	$\mathbf{u}_h^{n} \in G $ on $ S_K $ for all $K$ implies that $\mathbf{u}_h^{n+1}  \in G$ on $ S_K  $ for all $K.$
	\label{theorem:crkdg entropy boundedness 2d}
	\begin{proof}
The proof is analogous to the proof of Theorem \ref{theorem:crkdg entropy boundedness} and therefore we omit it.
	\end{proof}
\end{theorem}
\section{Numerical Experiments}
In this section, entropy-bounded cRKDG (EBcRKDG) is applied to one- and two- dimensional test cases to demonstrate its accuracy and ability to restrain spurious oscillations arising in discontinuous test cases. For two-dimensional problems we have chosen rectangular meshes. At each stage of the Rk method we have applied the positivity preserving limiter and the entropy-bound limiter consecutively after the stage update.
\subsection{One-dimensional test cases}
\begin{figure}[H]
\centering
\includegraphics[trim = 10 10 10 10, 
                width=\linewidth]{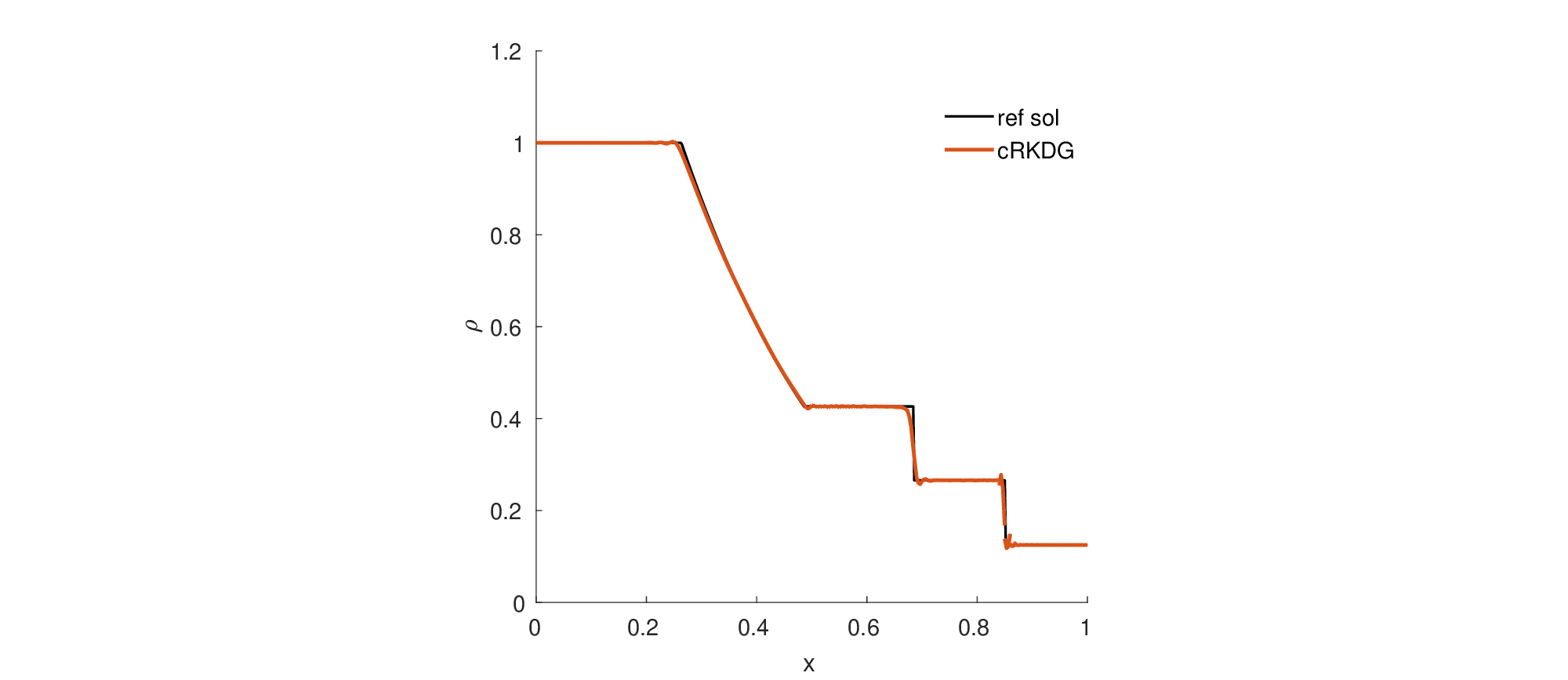}
 \caption{Plot of density for entropy-bounded cRKDG scheme, $k = 2$.}
 \label{figure: sod problem}
\end{figure}

\begin{example}
We consider the Sod's shock tube problem over the domain $x \in [0, 1] $ with following Riemann type initial conditions
\begin{align}
 (\rho, v, p)^T = \left\{ \begin{array}{l}
 	(1.0 \,, \, 0.0  \, , \,  1.0)^T \hspace{1.5cm} x \le 0.5, \\
 	(0.125, \, 0.0 , \, 0.1)^T \hspace{1.3cm} x > 0.5. 
 \end{array}                     \right.
\end{align} 
The domain is discretized using $100$ elements and the final time of computation is $ T = 0.2 .$ For limiting the solution the entropy bound is chosen to be $s_0 = \ln(1) $ i.e. the minimum entropy of the initial condition. In Figure \ref{figure: sod problem} we have shown the plot of entropy-bounded cRKDG scheme with $k = 2$ and a third order RK in Butcher form.  
\end{example}
\section{Summary and Conclusion}
In this paper we developed an entropy-bounded compact Runge-Kutta discontinuous Galerkin method(EBcRKDG) for hyperbolic conservation laws. It was first shown by Tadmor that the entropy solution of Euler equations of compressible gas dynamics satisfies a minimum-entropy principle \cite{Tadmor}. It was also proven by the author that the solutions of three point Finite volume/ Finite difference scheme 
like Lax-Friedrichs scheme and Godunov scheme also satisfy a similar entropy principle and its entropy is bounded below by the entropy of the solution in the domain of the dependence. An entropy-bounded DG method for Euler equations was proposed by Zhang et al.\cite{Zhang entropy} using a scaling limiter and Newton-Raphson iterative procedure. Later, Lv et al. proposed a more general EBDG scheme which avoids expensive iterative procedure, use local entropy bounds and could be extended to arbitrary elements in multi-dimension. However, the techniques developed by earlier researchers are only applicable to the DG schemes which are derived from the method of lines and use convexity arguments to extend the entropy-boundedness from forward-Euler time discretization schemes to one with strong-stability preserving (SSP) RK. Therefore, the cRKDG schemes, which are not derived from the standard method of lines and could only be expressed in Butcher form, deserve special attention in this regard. We observed that the element averaged solution with local derivative operator still satisfy an entropy principle similar to that with DG operator. The bound preserving framework of \cite{bound preserving crkdg} allowed us to express every stage update in RK method as a convex combinations of forward-Euler updates of previous stage values involving both local derivative operator and the DG operator. Finally, applying the entropy-bound limiter and the positivity preserving limiter furnish the entropy principle for the fully-discrete schemes.
		
		\clearpage
	\end{document}